\documentclass{article}
\usepackage{amssymb}
\usepackage[cp1251]{inputenc}
\usepackage[russian]{babel}
\usepackage{amsmath}
\usepackage{graphics}
\usepackage{epsfig}
\usepackage[all,knot]{xy}
\xyoption{arc}
\usepackage{color}

\def\mapr#1{\smash{\mathop{\buildrel{#1}\over\longrightarrow}}}

\def\cA{{\cal A}}
\def\cB{{\cal B}}
\def\cE{{\cal E}}
\def\rU{{\mathfrak U}}
\def\beq#1{\begin{equation}\label{#1}}
\def\eeq{\end{equation}}
\def\const{\mbox{\bf const }}
\newtheorem{theorem}{Тheorem}
\def\Aut{\hbox{\bf Aut}}

\title{Transitive Lie algebroids - categorical point of  view
}
\author{A.S.Mishchenko\\ Moscow State University \\ }
\date{}

\begin{document}
\maketitle

\section* {Introduction}

Transitive Lie algebroids have specific properties that allow to look at the
transitive Lie algebroid as an element of the object of a homotopy functor.
Roughly speaking each transitive Lie algebroids can be described as a vector
bundle over the tangent bundle of the manifold which is endowed with
additional structures. Therefore transitive Lie algebroids admits a
construction of inverse image generated by a smooth mapping of smooth
manifolds. The construction can be managed as a homotopy functor from the
category of smooth manifolds to the transitive Lie algebroids.     The
intention of this article is to make a classification of transitive Lie
algebroids and on this basis to construct a classifying space. The
realization of the intention allows to describe characteristic classes of
transitive Lie algebroids form the point of view a natural transformation of
functors similar to the classical abstract characteristic classes for vector
bundles.

\section{Definitions and formulation of the problem}

Given smooth manifold $ M $ let

$$
E \mapr{a} TM \mapr{p_{T}} M
$$
be a vector bundle over $ TM$ with fiber $ g$, $ p_{E}=p_{T}\cdot a$.
So we have a commutative diagram of two vector bundles
$$
\xymatrix{
E\ar[r]^{a}\ar[d]_{p_{E}}&TM\ar[d]^{p_{T}}\\
M\ar[r]&M
}
$$
The diagram is endowed with additional structure (commutator braces) and then is called
(\cite{Mck-05}, definition 3.3.1, \cite{Kub-91e}, definition 1.1.1)
transitive Lie algebroid
$$
\cA = \left\{
\makebox(70,25)[t]{\xymatrix{
E\ar[r]^{a}\ar[d]_{p_{E}}&TM\ar[d]^{p_{T}}\\
M\ar[r]&M
}}; \{\bullet,\bullet\}
\right\}.
$$

Let $ f:M'\mapr{}M $ be a smooth map. Then one can define an inverse image
(pullback) of the Lie algebroid (\cite{Mck-05}, page 156, \cite{Kub-91e}, definition 1.1.4),
$ f^{!!}(\cA)$. This means that given the finite dimensional Lie algebra $g$
there is the functor $\cA$ such that
with any manifold $M$ it assigns the family $\cA(M)$ of all transitive  Lie algebroids
with fixed Lie algebra $g$.

In the dissertation \cite{Wa-2007} the following statement was proved:
Each transitive Lie algebroid is trivial, that is there is a trivialization
of vector bundles $E, TM, \ker a=\bar{g}$ such that
$$
E\approx TM\oplus\bar{g},
$$
and the Lie bracket is defined by the formula:
$$
[(X,u),(Y,v)]= ([X,Y], [u,v]+X(v)-Y(u)).
$$

Then using the construction of pullback and the idea by Allen Hatcher
\cite{Ha-2003} one can prove that the functor $\cA$ is homotopic functor.
More exactly for two homotopic smooth maps $f_{0},f_{1}:M_{1}\mapr{}M_{2}$
and for the transitive  Lie algebroid
$$\left(E\mapr{a}TM_{2}\mapr{}M_{2};\{\bullet,\bullet\}\right)$$
two inverse images $f_{0}^{!!}(\cE)$ and $f_{1}^{!!}(\cA)$ are isomorphic.

Hence there is a final classifying space $\cB_{g}$ such that
the family of all transitive  Lie algebroids with fixed Lie algebra $g$ over the manifold $M$
has one-to-one correspondence with the family of homotopy classes of continuous maps $[M,\cB_{g}]$:
$$
\cA(M)\approx [M,\cB_{g}].
$$

Using this observation one can describe the family of all characteristic classes of
a transitive  Lie algebroids in terms of cohomologies of the classifying space $\cB_{g}$. Really,
from the point of view of category theory a characteristic class $\alpha$ is a natural transformation
from the functor $\cA$ to the cohomology functor $H^{*}$. This means that for the transitive  Lie algebroid
$\cE=\left(E\mapr{a}TM\mapr{}M;\{\bullet,\bullet\}\right)$ the value of the characteristic class
$\alpha(\cE)$ is a cohomology class
$$
\alpha(\cE)\in H^{*}(M),
$$
such that for smooth map $f:M_{1}\mapr{}M$ we have
$$
\alpha(f_{0}^{!!}(\cE))=f^{*}(\alpha(\cE))\in H^{*}(M_{1}).
$$

Hence the family of all characteristic classes $\{\alpha\}$ for  transitive
Lie algebroids with fixed Lie
algebra $g$ has a one-to-one correspondence
with the cohomology group $H^{*}(\cB_{g})$.

On the base of these abstract considerations a natural problem can be
formulated.

{\bf Problem.} Given finite dimensional Lie algebra $g$ describe the
classifying space $\cB_{g}$
for  transitive Lie algebroids
in more or less understandable terms.

Below we suggest a way of solution the problem and consider some trivial
examples.

\section{Description of transitive Lie algebroids using transition functions}

Consider the trivial transitive Lie algebroids
$$
E\approx TM\oplus\bar{g}, \quad \bar{g}\approx M\times g,
$$
and the Lie bracket is defined by the formula:
$$
[(X,u),(Y,v)]= ([X,Y], [u,v]+X(v)-Y(u)),
$$
where $X,Y\in \Gamma^{\infty}(TM)$ are smooth vector fields,
$u,v\in \Gamma^{\infty}\bar{g}$ are smooth sections which are represented as
smooth vector functions with values in the Lie algebra $g$. Consider
a fiberwise isomorphism $\cA: E\mapr{}E$ that is identical on the second summands
and generates the Lie algebra homomorphism
$\cA:\Gamma^{\infty}(E)\mapr{}\Gamma^{\infty}(E)$. The isomorphism $\cA$ can
be written by formula:
$$\begin{array}{l}
(v,Y)=\cA(u,X); \\
(v,X)=(\varphi(x)(u(x))+\omega(X), X),
\end{array}
$$
where $\varphi(x):g\mapr{}g$ is a fiberwise map of the bundle $\bar{g}$,
and $\omega$ is a differential form with values in $g$.
The isomorphism $\cA$ can be expressed as a matrix
$$
\left(
\begin{array}{c}
v(x) \\Y
\end{array}
\right)=
\left(
\begin{array}{cc}
\varphi(x) & \omega\\
0 & 1
\end{array}
\right)\cdot\left(
\begin{array}{c}
u(x) \\ X
\end{array}
\right)
$$

From the property of that $\cA$ is a Lie algebra homomorphism:
$$
\cA([(X,u),(Y,v)])=[\cA(X,u),\cA(Y,v)]
$$
one has that
\beq{1}
\begin{array}{l}
\varphi(x)([u_{1}(x),u_{2}(x)])=[\varphi(x)(u_{1}(x)),\varphi(x)(u_{2}(x))]),\\\\
d\omega(X_{1},X_{2})+[\omega(X_{1}),\omega(X_{2})]=0,\\\\
d\varphi(X)(u)=[\varphi(u),\omega(X)].
\end{array}
\eeq

Consider an atlas of charts on the manifold $M$, $\{\rU_{\alpha}\}$,
$\bigcup\limits_{\alpha}U_{\alpha}=M$, and the trivializations
$E_{\alpha}\stackrel{\Phi_{\alpha}}{\approx} TU_{\alpha}\otimes(U_{\alpha}\times g)$
of the Lie algebroid $E$ over each chart $U_{\alpha}$ with the Lie brackets
defined by the formula
$$
[(X,u),(Y,v)]= ([X,Y], [u,v]+X(v)-Y(u)),
$$
for $X,Y\in \Gamma^{\infty}(TU_{\alpha})$,
$u,v\in \Gamma^{\infty}(U_{\alpha}\times g)$.

On the intersection of two charts
$U_{\alpha\beta}=U_{\alpha}\cap U_{\beta}$
we have the transition function
$$
\Phi_{\beta\alpha}=\Phi_{\beta}\Phi^{-1}_{\alpha}:
TU_{\alpha\beta}\otimes(U_{\alpha\beta}\times g)\mapr{}
TU_{\alpha\beta}\otimes(U_{\alpha\beta}\times g)
$$
which have the matrix form

$$
\left(
\begin{array}{c}
v(x) \\Y
\end{array}
\right)=
\Phi_{\beta\alpha}\left(\begin{array}{c}
u(x) \\ X
\end{array}\right)=
\left(
\begin{array}{cc}
\varphi_{\beta\alpha}(x) & \omega_{\beta\alpha}\\
0 & 1
\end{array}
\right)\cdot\left(
\begin{array}{c}
u(x) \\ X
\end{array}
\right).
$$

For another choice of trivializations $\Phi'_{\alpha}$ the correspondent transition
functions $\Phi'_{\beta\alpha}$ satisfy the homology condition:
$$
\Phi'_{\beta\alpha}=
{\rm H}_{\beta}\cdot\Phi_{\beta\alpha}\cdot{\rm H}^{-1}_{\alpha}
$$
$$
\begin{array}{l}
\left(
\begin{array}{cc}
\varphi'_{\beta\alpha}(x) & \omega'_{\beta\alpha}\\
0 & 1
\end{array}
\right)= \\\\=
\left(
\begin{array}{cc}
\eta_{\beta}(x) & \mu_{\beta}\\
0 & 1
\end{array}
\right)\cdot
\left(
\begin{array}{cc}
\varphi_{\beta\alpha}(x) & \omega_{\beta\alpha}\\
0 & 1
\end{array}
\right)\cdot
\left(
\begin{array}{cc}
\eta^{-1}_{\alpha}(x) & -\eta^{-1}_{\alpha}\mu_{\alpha}\\
0 & 1
\end{array}
\right)
\end{array}
$$

or

$$
\begin{array}{l}
\left(
\begin{array}{cc}
\varphi'_{\beta\alpha}(x) & \omega'_{\beta\alpha}\\
0 & 1
\end{array}
\right)= \\\\ =
\left(
\begin{array}{cc}
\eta_{\beta}(x)\varphi_{\beta\alpha}(x)\eta^{-1}_{\alpha}(x) &
-\eta_{\beta}(x)\varphi_{\beta\alpha}(x)\eta^{-1}_{\alpha}(x)\mu_{\alpha}+
\eta_{\beta}(x)\omega_{\beta\alpha}+\mu_{\beta}
\\
0 & 1
\end{array}
\right),
\end{array}
$$

or
$$
\begin{array}{l}
\varphi'_{\beta\alpha}(x)= \eta_{\beta}(x)\varphi_{\beta\alpha}(x)\eta^{-1}_{\alpha}(x),\\\\
\omega'_{\beta\alpha}=-\eta_{\beta}(x)\varphi_{\beta\alpha}(x)\eta^{-1}_{\alpha}(x)\mu_{\alpha}+
\eta_{\beta}(x)\omega_{\beta\alpha}+\mu_{\beta}.
\end{array}
$$

The elements $\eta_{\beta}$ and $\mu_{\beta}$ satisfy similar (\ref{1}) conditions:

$$
\begin{array}{l}
\eta_{\beta}(x)([u_{1}(x),u_{2}(x)])=[\eta_{\beta}(x)(u_{1}(x)),\eta_{\beta}(x)(u_{2}(x))]),\\\\
d\mu_{\beta}(X_{1},X_{2})+[\mu_{\beta}(X_{1}),\mu_{\beta}(X_{2})]=0,\\\\
d\eta_{\beta}(X)(u)=[\eta_{\beta}(u),\mu_{\beta}(X)].
\end{array}
$$

\section{Case of commutative Lie algebra $g$ }

In commutative case the conditions (\ref{1}) have for simple form:
\beq{2}
\begin{array}{l}
\varphi_{\beta\alpha}(x)([u_{1}(x),u_{2}(x)])=
[\varphi_{\beta\alpha}(x)(u_{1}(x)),\varphi_{\beta\alpha}(x)(u_{2}(x))]),\\\\
d\omega_{\beta\alpha}(X_{1},X_{2})=0,\\\\
d\varphi_{\beta\alpha}(X)(u)=0.
\end{array}
\eeq
Hence
$$
\varphi_{\beta\alpha}(x)=\const.
$$

This means that the vector bundle $\bar{g}$ is flat and
the family $\omega=\{\omega_{\beta\alpha}\}$ defines a Cech cochain
$$
\omega\in C^{1}(\rU,\Omega^{1}(\bar{g}))
$$
in the bigraded Cech complex
$$C^{*,*}=\left\{\bigoplus C^{i}(\rU, \Omega^{j}(\bar{g}); d=d'+d''\right\}$$
where $\rU=\{U_{\alpha}\}$ is the atlas of charts.

One has
$$
d'(\omega)=0;\quad d''(\omega)=0.
$$
Hence $\omega$ defines cohomology class
$$
[\omega]\in H^{2}(M;\bar{g}).
$$

Therefore we have the following
\begin{theorem} The classification of all transitive Lie algebroids
with fixed commutative
Lie algebra $g$ over the manifold $M$
is determined by a flat Lie algebra bundle $\bar{g}$ over $M$ and
a 2-dimensional cohomology class $[\omega]\in H^{2}(M;\bar{g})$.
\end{theorem}

\section{Some general properties}

In common case we can say that a little bit about the transition functions on the
level of homology groups $H_{*}(g)$ of the Lie algebra $g$. Since each transition
function $\varphi_{\beta\alpha}(x)$ is the homomorphism of the Lie algebra $g$,
that is $\varphi_{\beta\alpha}(x)\in \Aut(g)$, the cocycle
$\{\varphi_{\beta\alpha}(x)\}$ generate associated bundles with fibers $H_{*}(g)$,
say, $\overline{H_{*}(g)}$, and bundles with fibers $H^{*}(g)$,$\overline{H^{*}(g)}$.
The properties (\ref{1}) imply that all bundles
$\overline{H_{*}(g)}$ and $\overline{H^{*}(g)}$
are flat. In particular the differential forms
$\omega_{\beta\alpha}\in \Omega^{1}(U_{\alpha\beta};\bar{g})$ generates the cocycle
$$
\overline{\omega}=
\{\overline{\omega}_{\beta\alpha}\}\in C^{1}(\rU,\overline{H_{1}(g)})=
\bigoplus_{\alpha\beta}\Omega^{1}(U_{\alpha\beta};\overline{H_{1}(g)}),
$$
that is
$$
d'(\overline{\omega})=0,
$$
$$
d''(\overline{\omega})=0.
$$
This means that the cocycle $\overline{\omega}$ induces a cohomology
class
$$
[\overline{\omega}]\in H^{2}\left(M;\overline{H_{1}(g)}\right).
$$

The foregoing consideration creates a conjecture that classification of the transitive
Lie algebroid $E$ induces by two things: the Lie algebra bundle with structural group
$\widetilde{\Aut(g)}$ with special topology and the cohomology class
$[\overline{\omega}]\in H^{2}\left(M;\overline{H_{1}(g)}\right)$. The special topology
in the group $\Aut(g)$ is defined as a minimal topology, which is more fine topology
than the classical
topology in $\Aut(g)$ and such that all homomorphisms
$$
\Aut(g)\mapr{}\Aut(H_{k}(g))_{discrete}
$$
are continuous.

This paper is partly supported  by grants RBRF No
05-01-00923-a,07-01-91555-NNIO-a, 10-01-92601-KO-a and project No. RNP.2.1.1.5055

\end{document}